\documentclass[twocolumn]{revtex4-1}
\pdfoutput=1

\usepackage{amsmath}
\usepackage{amssymb}
\usepackage{amsfonts}
\usepackage{graphicx}

\begin{document}

\title{\Large A New Equilateral Triangle? }
\author{Martin Buysse}
\email[]{martin.buysse@uclouvain.be}
\affiliation{{Facult\'e d'architecture, d'ing\'enierie architecturale, d'urbanisme -- LOCI}, {UCLouvain}}

\begin{abstract} 
\noindent In any triangle, the perpendicular side bisectors meet the corresponding internal angle bisectors on the circumcircle. If we take those three points as the vertices of a new triangle and repeat the operation indefinitly, we end up in the limit with a par of equilateral triangles whose sides are parallel to the sides of the Morley triangle of the initial triangle. 
\end{abstract}

\maketitle

\noindent Of all the triangles, the equilateral is the most symmetrical. Scalene triangles, on the other hand, have no symmetry. This is probably why the univocal association of an equilateral triangle, having its determined size, localization and orientation, to each scalene triangle -- actually to any reference triangle --, is so fascinating. 
The triangle of Morley, whose vertices are the three intersection points of adjacent internal angle trisectors, is a striking example \cite{Morley:1924,Oakley:1978}. Dozen were actually found out by considering all the trisectors and their intersections. Still, there are others. Two of them are attributed, albeit probably abusively, to Napoleon Bonaparte \cite{Rutherford:1825,Grunbaum:2012}! The theorem states that if equilateral triangles are erected on the sides of any triangle, the centers of those three triangles themselves form an equilateral triangle: the outer Napoleon triangle if the triangles are erected outwards, and the inner one if they are erected inwards. 
\begin{figure}[h]
  \centering
  \includegraphics[width=86mm,keepaspectratio,trim=0 90 0 40,clip=true]{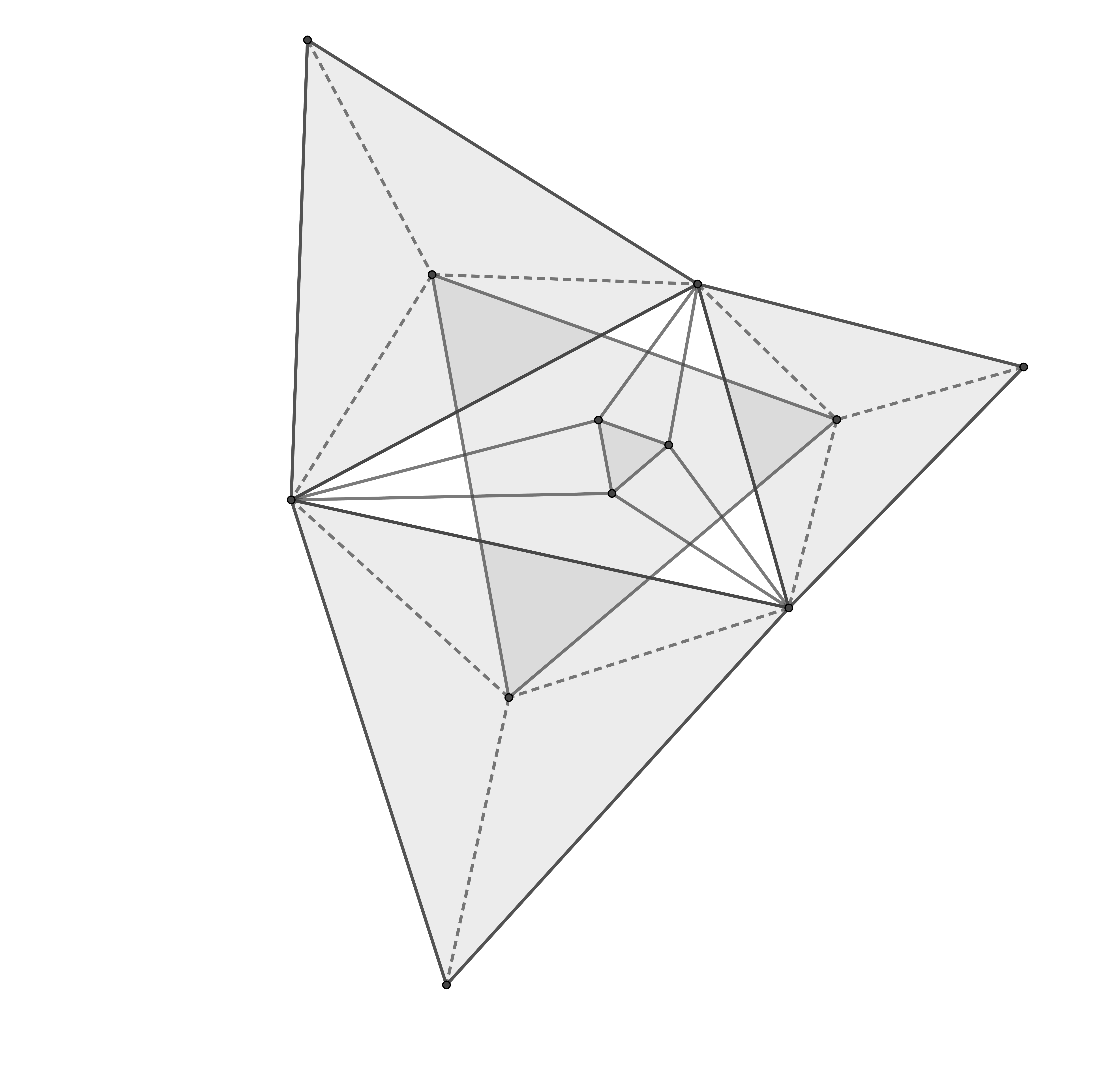}%
  \caption{Morley \& outer Napoleon triangles}
\end{figure}

There are some pedal, antipedal, circumcevian triangles, and many more -- see \cite{McCartin:2012} for a review or \cite{Oai:2015,*Oai:2018a,*Oai:2018b,Dasari:2019} for recent attempts. We found another one -- actually two of them -- , that might have been unknown until now, even if it has two known cousins: one is smaller, infinitely; and the other one larger, infinitely too. We will come to them later. 

Let $a$, $b$, $c$ be the sidelengths of a reference triangle; $\alpha$, $\beta$, $\gamma$ its opposite angles; and $\ell_a$, $\ell_b$, $\ell_c$, the lengths of the corresponding arcs on the circumcircle. The perpendicular side bisectors also bisect the corresponding arcs; so do the opposite angle bisectors. We consider the triangle whose vertices are those three intersection points, its three perpendicular and internal angle bisectors and their intersection points on the circumcircle, that we use to build a new triangle, and so on. The sequence of triangles converges to a par of equilateral triangles symmetrical with respect to the circumcenter: the triangles of even $n$ rank converge to the first equilateral triangle, and the triangles of odd rank to the second one. 
\begin{figure}[h]
  \centering
  \includegraphics[width=43mm,keepaspectratio]{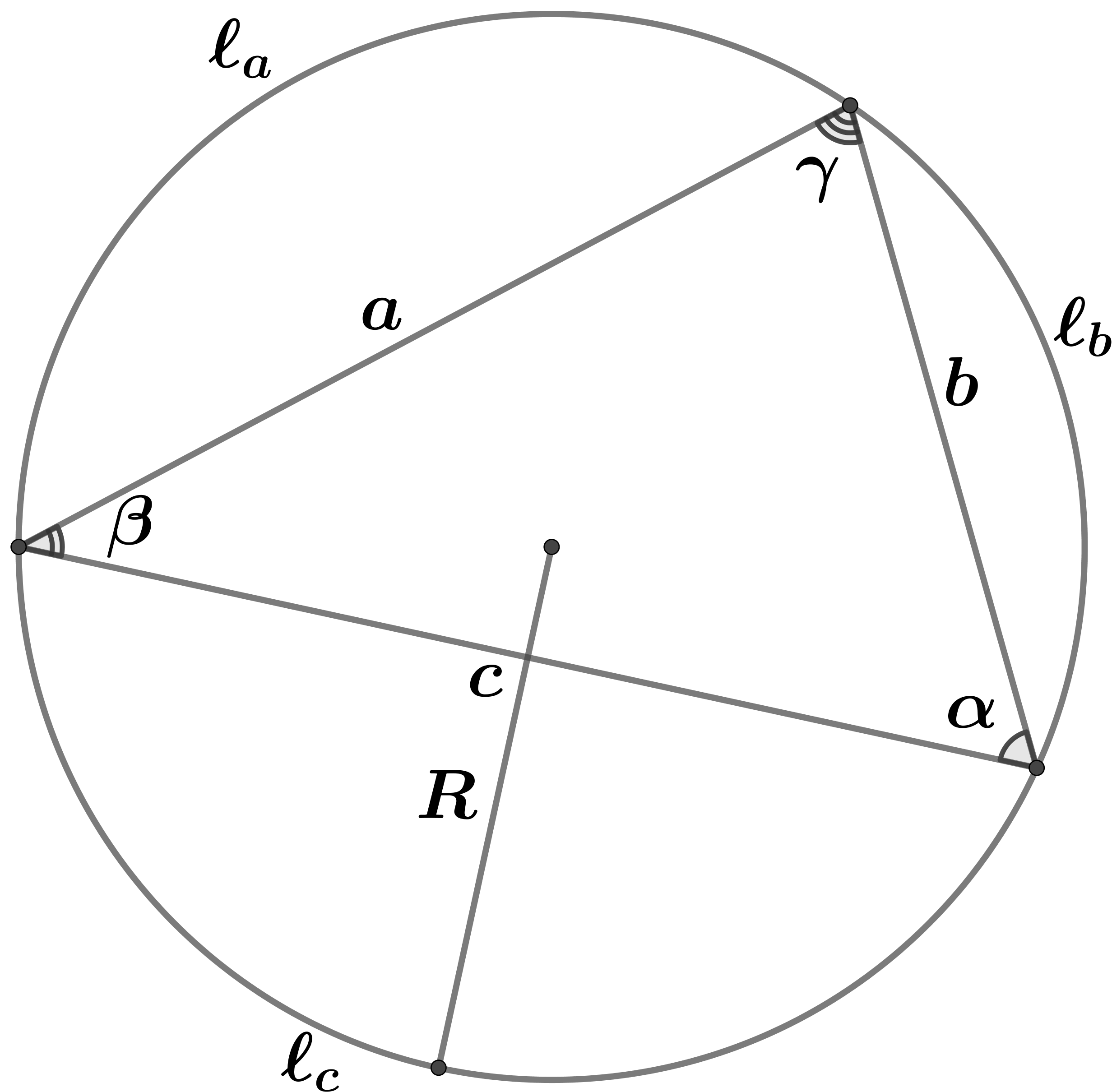}%
  \includegraphics[width=43mm,keepaspectratio]{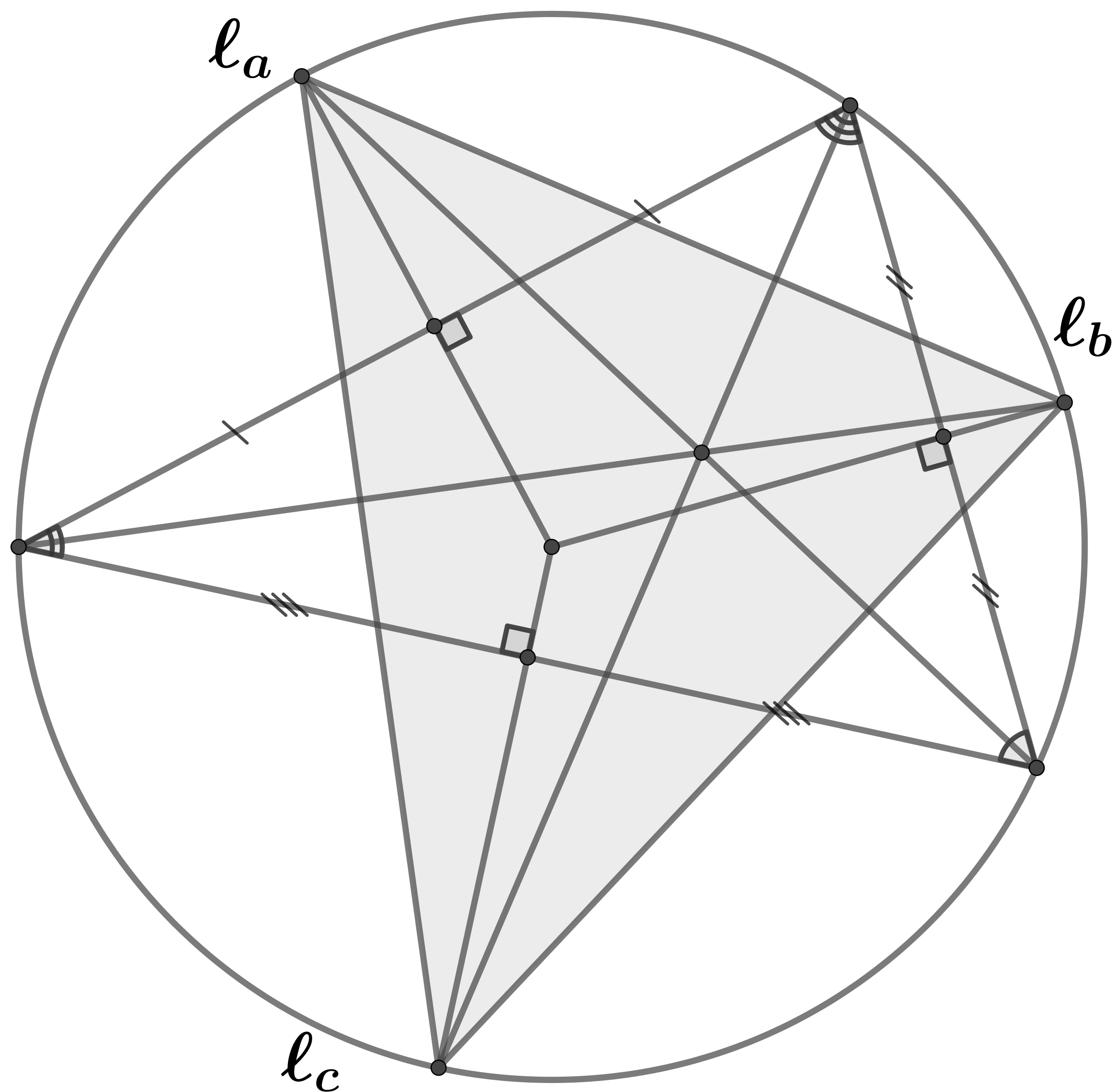}
  \includegraphics[width=43mm,keepaspectratio]{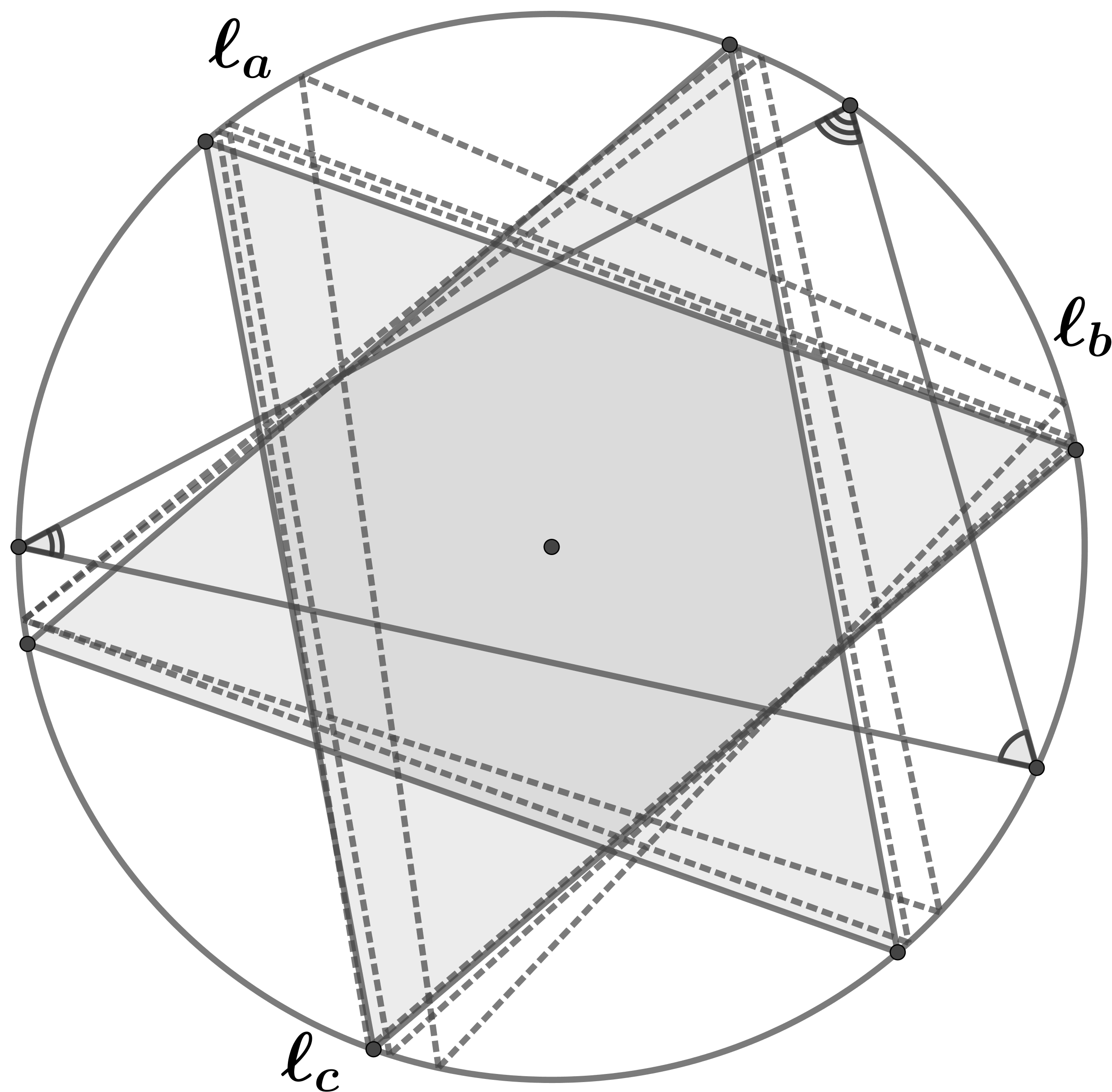}%
  \includegraphics[width=43mm,keepaspectratio]{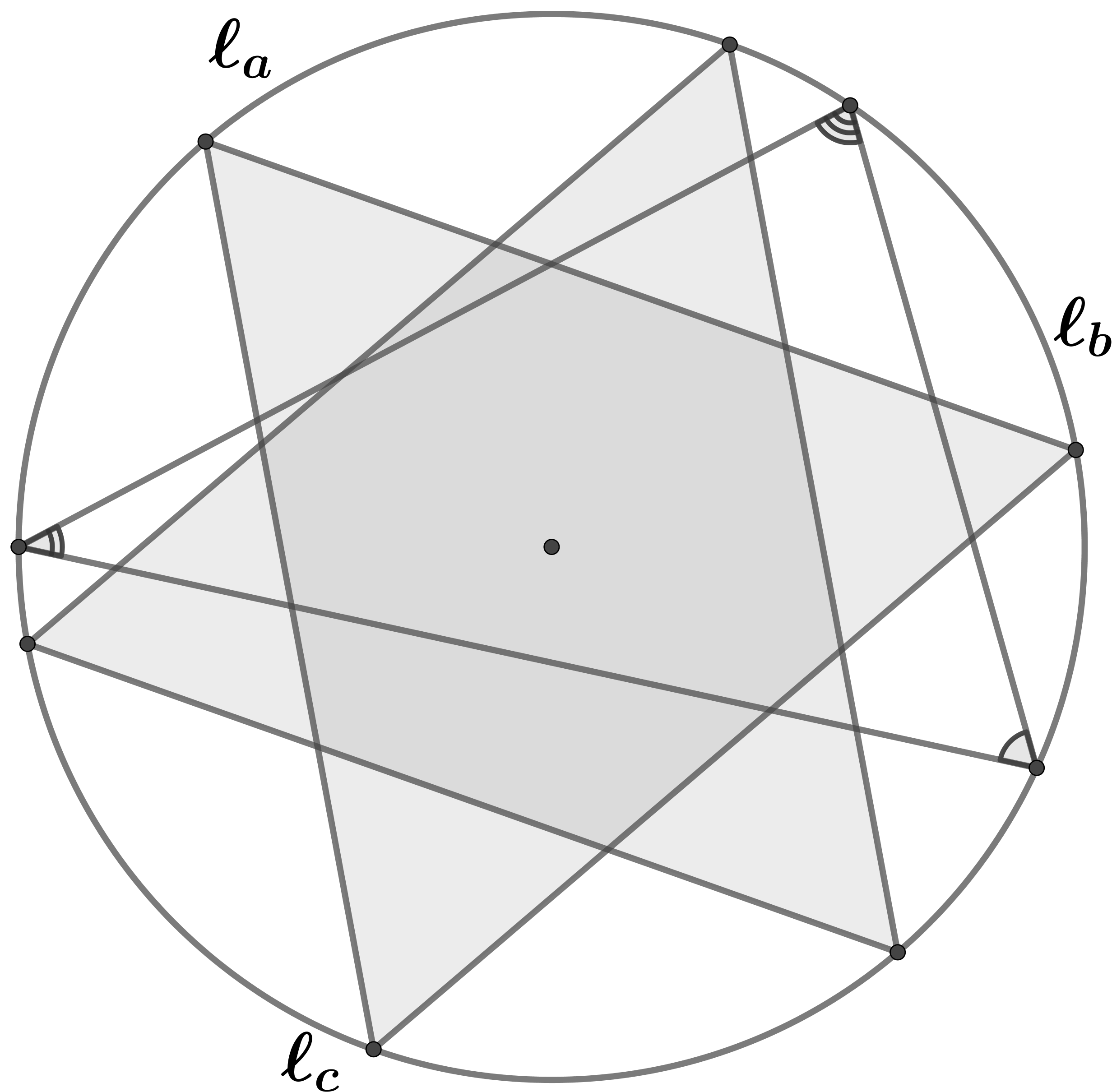}%
  \caption{Iteration process from a reference triangle}
\end{figure}
This result can be proven by computing the arclengths of the triangles of this infinite sequence: 
\begin{align}
&\;(\ell_a,\ell_b,\ell_c)&n=0
\nonumber\\
\rightarrow &\left(\frac{\ell_b+\ell_c}{2},\frac{\ell_a+\ell_c}{2},\frac{\ell_a+\ell_b}{2}\right)&n=1
\nonumber\\
\rightarrow &\left(\frac{2\ell_a+\ell_b+\ell_c}{4},\frac{\ell_a+2\ell_b+\ell_c}{4},\frac{\ell_a+\ell_b+2\ell_c}{4}\right)&\vdots\hspace{.4cm}
\nonumber
\\
\rightarrow &\left(\frac{2\ell_a+3\ell_b+3\ell_c}{8},\frac{3\ell_a+2\ell_b+3\ell_c}{8},\frac{3\ell_a+3\ell_b+2\ell_c}{8}\right)\hspace{-1cm}&
\nonumber
\end{align}

\noindent that is, at the $n$th iteration
\begin{align}
&
\left(
\frac{\frac{2^n+2}{3}\ell_a+\frac{2^n-1}{3}\ell_b+\frac{2^n-1}{3}\ell_c}{2^n},
\ldots,\ldots
\right)
\nonumber\\
&
\left(
\frac{\frac{2^n-2}{3}\ell_a+\frac{2^n+1}{3}\ell_b+\frac{2^n+1}{3}\ell_c}{2^n},
\ldots,\ldots
\right)
\nonumber
\end{align}

\noindent respectively for the triangles of even and odd $n$ rank. When $n$ gets arbitrarily large, both of them tend to 
\begin{equation*}
\left(
\frac{\ell_a+\ell_b+\ell_c}{3},
\frac{\ell_a+\ell_b+\ell_c}{3},
\frac{\ell_a+\ell_b+\ell_c}{3}
\right)
\end{equation*}

\noindent i.e. the dimensions of an equilateral triangle. 

\begin{figure}[h]
  \centering
  \includegraphics[width=86mm,keepaspectratio,trim=0 0 0 30,clip=true]{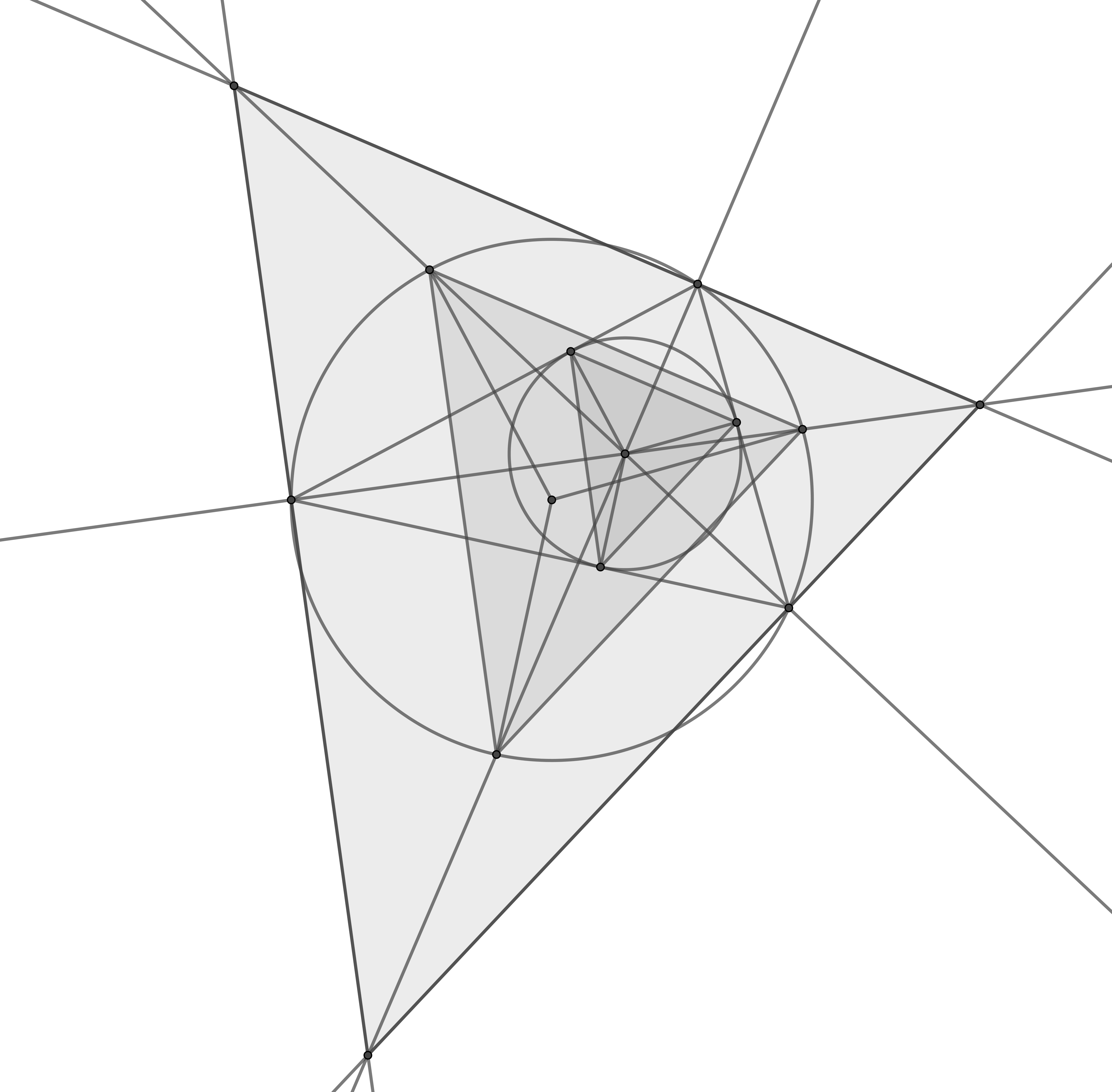}%
  \caption{Contact, excentral \& new triangle: first iteration}
\end{figure}

To indroduce the cousins, recall that the excentral triangle of a reference triangle is the triangle whose vertices are the excenters, that is the intersection points of the external angle bisectors. Since the external angle bisectors are perpendicular to the internal ones, the excentral triangle is similar to our first iteration ($n=1$) triangle, the excentral triangle of the excentral triangle to our second triangle, etc. It is known that this sequence of excentral triangles approaches an equilateral triangle \cite{Johnson:2013,Villiers2014}. 

The contact triangle, whose vertices are the contact points of the reference triangle with its incircle, is also similar to our first iteration triangle, the contact triangle of the contact triangle to our second triangle, etc. This sequence of triangles approaching equiangularity is discussed in \cite{Chang:1997}, \cite{Villiers2014}. The triangles of the three sequences have, at each step, their sides parallel to each other.  

\vspace{3mm}

The orientation of the equilateral triangles has yet to be determined. In the sequence, we focus on the triangles of even rank and calculate the length $\ell_{ab}$ of the circumcircle arc between one vertex, say the $(a,b)$ vertex with $a\geq b$, of the reference triangle, and the corresponding vertex of the rank $n$ triangle. This arclength is obviously $0$ when $n=0$. For $n=2$, $n=4$ and $n\rightarrow\infty$, we have
\begin{align}
\ell_{ab,n=2}&=\frac{\ell_a-\ell_b}{4}
\qquad\qquad
\ell_{ab,n=4}=\frac{\ell_a-\ell_b}{4}+\frac{\ell_a-\ell_b}{16}
\nonumber\\
\ell_{ab,n\rightarrow\infty}&=\frac{\ell_a-\ell_b}{4}\sum_{k=0}^{\infty}\left(\frac{1}{4}\right)^k=\frac{\ell_a-\ell_b}{4}\frac{1}{1-\frac{1}{4}}=\frac{\ell_a-\ell_b}{3}
\nonumber
\end{align}

\noindent Since $2R=\ell_a/\alpha=\ell_b/\beta=\ell_c/\gamma$, where $R$ is the length of the circumradius, this result means that the angle between the corresponding sides of the reference and the final even triangle is given by $(\alpha-\beta)/3$, that is, the same orientation as the Morley and outer Napoleon triangles. 
\begin{figure}[h]
  \centering
  \includegraphics[width=43mm,keepaspectratio]{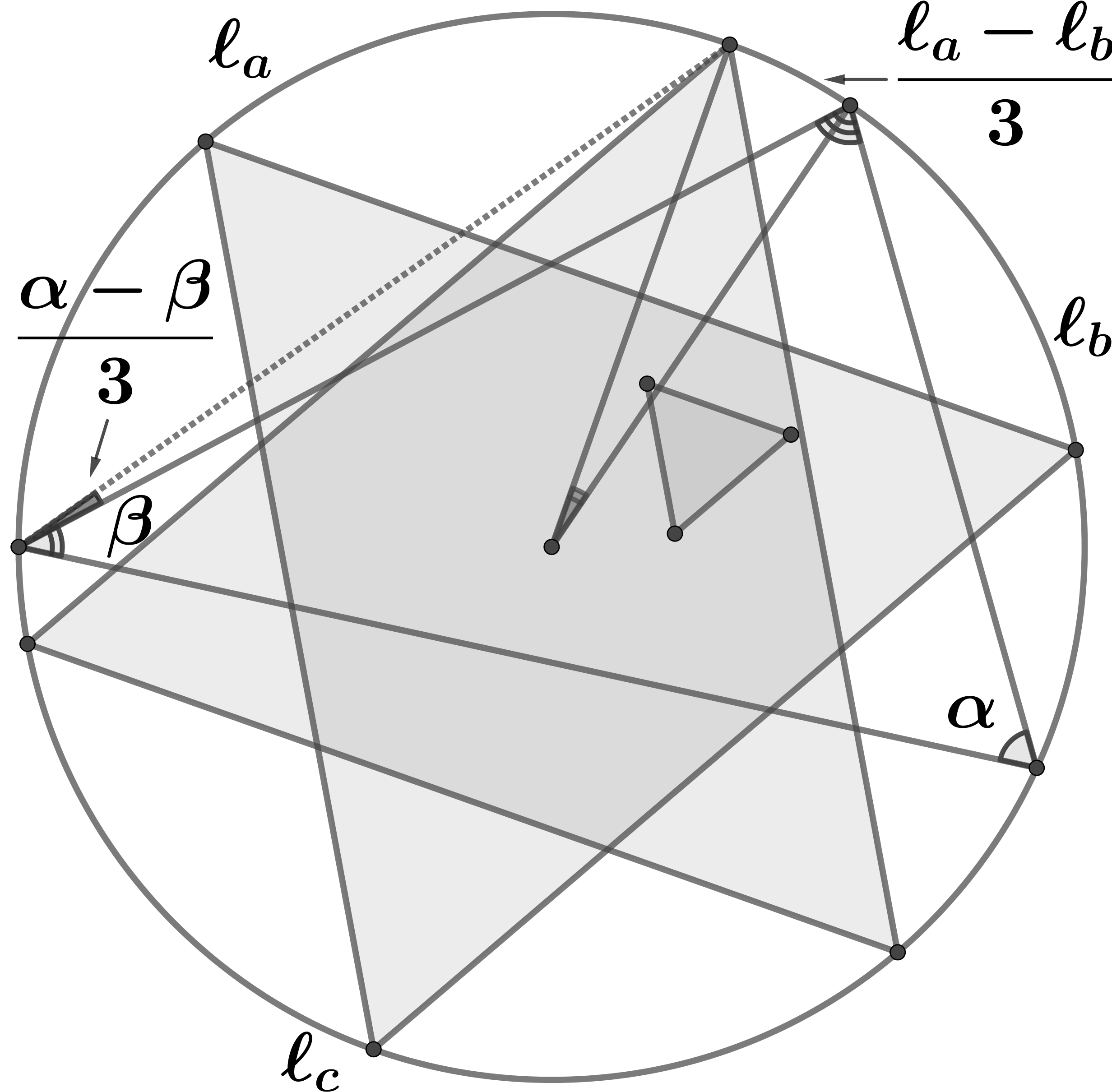}%
  \includegraphics[width=43mm,keepaspectratio]{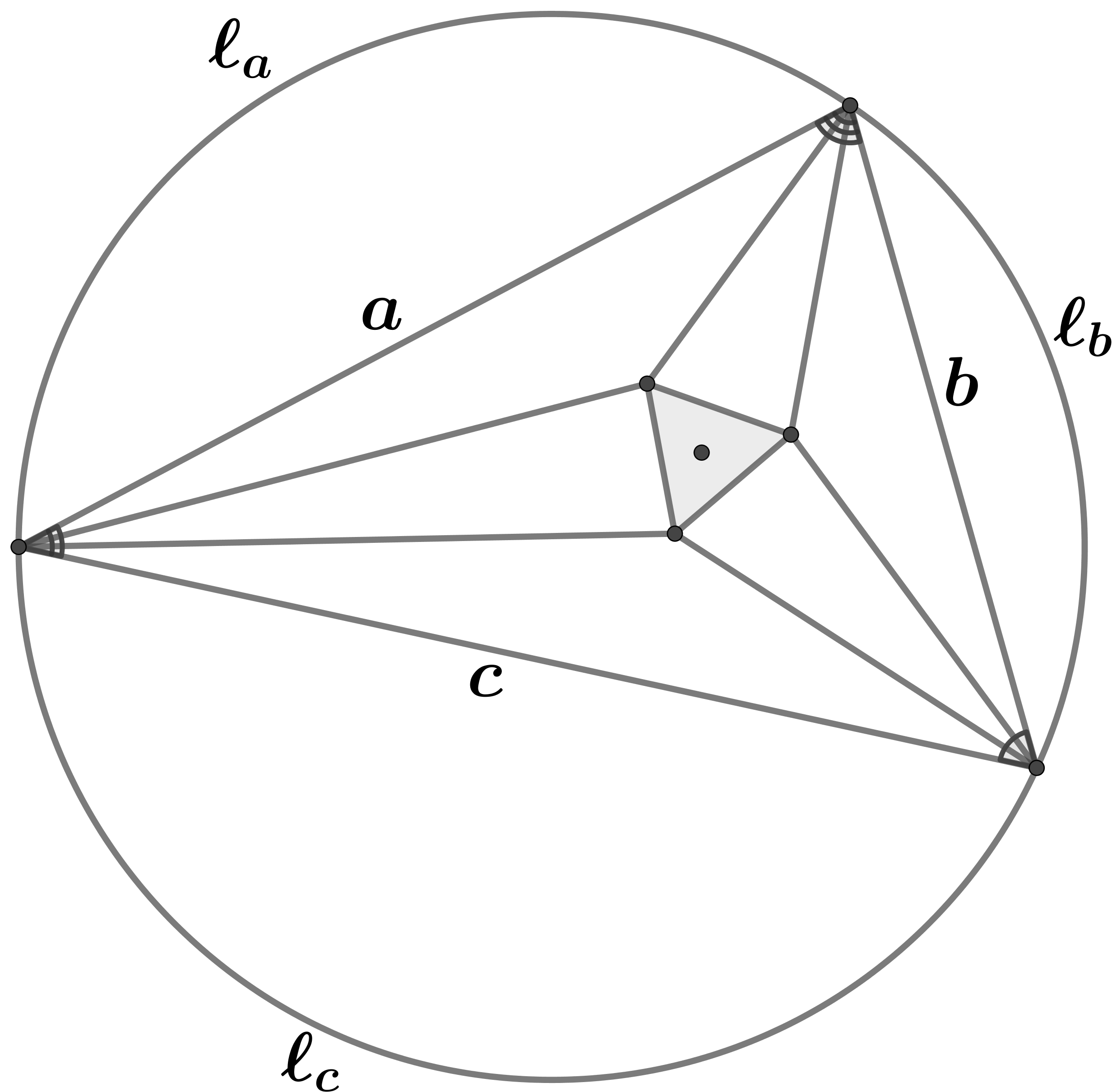}%
  \caption{Orientation of the new triangles \& Morley triangle}
\end{figure}
 
Its sidelength is equal to $\sqrt{3}R$, which can be expressed in terms of $a$, $b$ and $c$ since the area of the reference triangle is given by $a b\sin\gamma /2=abc/(4R)$ on the one hand, and by the Heron formula on the other hand. 

\bibliographystyle{prsty}
\bibliography{BibliMath}

\end{document}